 
\def\L{{\cal L}}
\baselineskip=14pt
\parskip=10pt
\def\Tilde{\char126\relax}
\def\halmos{\hbox{\vrule height0.15cm width0.01cm\vbox{\hrule height
 0.01cm width0.2cm \vskip0.15cm \hrule height 0.01cm width0.2cm}\vrule
 height0.15cm width 0.01cm}}
\font\eightrm=cmr8  
\font\eighttt=cmtt8
\magnification=\magstephalf

\parindent=0pt
\overfullrule=0in
\bf
\centerline{The  ABSTRACT LACE EXPANSION }
\it
\centerline{Doron ZEILBERGER \footnote{$^1$}
{\baselineskip=9pt
\eightrm  \raggedright
Department of Mathematics, Temple University,
Philadelphia, PA 19122, USA. 
{\baselineskip=9pt
\eighttt zeilberg@math.temple.edu ; \break http://www.math.temple.edu/\Tilde 
zeilberg. }
Supported in part by the NSF. Dec. 5, 1996.
}}
\medskip
\qquad\qquad\qquad\qquad\qquad\qquad\qquad
\qquad\qquad\qquad\qquad\qquad\qquad
{\it For Erd\H os P\'al, In Memoriam}
\medskip
\rm
{\bf Abstract:}
David Brydges and Thomas Spencer's Lace Expansion
is abstracted, and it is shown how it sometimes gives rise to sieves.
\bigskip
 
{\bf LACES}
 
{\bf Definition:} Let $P$ be a finite set of {\it properties}.
A mapping $l$ that assigns to any subset $S \subset P$ another
subset $l(S)$, is called a  {\it lace-map}, if for all 
$S,G,S_1,S_2 \subset P$: \hfill\break 
$ (i)\,\,l(S) \subset S \,\,; \quad
(ii)\,\,l(S) \subset G \subset S \Rightarrow l(G)=l(S)\,\,; \quad
(iii)\,\, l(S_1)=l(S_2) \Rightarrow l(S_1 \cup S_2)=l(S_1) \,\,.$
 
A set $L$ for which $l(L)=L$ is called a {\it lace}.
By applying $(ii)$ to $G=l(S)$, it is seen that $l(l(S))=l(S)$,
for any set of properties $S$, hence $l(S)$ is always a lace,
and $l$ is a projection: $l^2=l$.
 
If $L$ is a lace then, by $(iii)$,
there exists a set $C(L) \subset P\backslash L$ such
that 
$$
\{ S \subset P \,\,|\,\, l(S)=L \}=\{ S
\,\,|\,\, L \subset S \subset L \cup C(L)\}
\quad .
$$
The
set $C(L)$ is called the {\it set of  properties
compatible with $L$}. The collection of laces will be denoted by $\L$.
For any lace $L$, obviously
$C(L)= \{ p \in P \backslash L \,\,|\,\, l(L \cup \{ p \} )=L \}.$
 
{\bf Theorem}: Let $X$ be a set of elements each of which
possesses a subset of the properties of $P$. Let
$wt$ be any function on $X$ (in particular the counting function
$wt(x)\equiv 1$). For any lace $L$ define
$N(L)$ to be the sum of the weights of those elements of $X$
that {\it definitely have} all the properties of $L$ and
{\it definitely don't have} any of the properties in $C(L)$.
Then the sum of the weights of those elements of $X$ that
have none of the properties of $P$, $N_0(X)$, is given by:
$$
N_0(X)=
\sum_{L \in \L} (-1)^{|L|} N(L) \quad .
\eqno(Lace\_Expansion)
$$
 
{\bf Proof:} For each property $p \in P$, assign a variable $Y_p$.
Since every subset $S$ of $P$ has a unique lace $L=l(S)$, and
by $(iii)$, the collection of subsets $G$ for which $l(G)=L$ 
consist of the interval (in the Boolean lattice) between
$L$ and $L \cup C(L)$, we have
$$
\prod_{p \in P} (1+ Y_p)=
\sum_{S \subset P} \prod_{s \in S} Y_s =
\sum_{L \in \L} \sum_{{{S;}\atop { l(S)=L}}} \prod_{s \in S} Y_s =
$$
$$
\sum_{L \in L} \prod_{s \in L} Y_s 
\sum_{{{S;}\atop { l(S)=L}}} \prod_{s \in S\backslash L} Y_s =
\sum_{L \in \L} \prod_{s \in L} Y_s
\prod_{s \in C(L)} (1+ Y_s) \quad .
\eqno(*)
$$
 
For every $x \in X$ and $p \in P$, let
$\chi_p(x)=-1$ if $x$ has property $p$ and $0$ if $x$ doesn't.
Since $0^r=0$ when $r>0$, while $0^0=1$, we have 
$$
N_0(X)=
\sum_{n X} wt(x) \prod_{ p \in P} (1+ \chi_p(x))=
\sum_{x \in X} wt(x) \sum_{L \in \L} \prod_{p \in L} \chi_p(x) 
\prod_{p \in C(L)} (1+ \chi_p(x))
=
$$
$$
\sum_{L \in \L} \sum_{x \in X} wt(x) \prod_{p \in L} \chi_p(x)
\prod_{p \in C(L)} (1+ \chi_p(x))=
\sum_{L \in \L} (-1)^{|L|} N(L) \quad \halmos \,\,.
$$

{\bf SIEVES} 
 
Let's call a lace {\it saturated} if its set of compatible properties,
$C(L)$, is empty. In this case $N(L)$ is simply the sum of the
weights of the elements of $X$ that definitely 
have all the properties in $L$
(and possibly others). Let $\L_s$ be the collection of
saturated laces.
If $l$ is such that the parity of the cardinalities of all
unsaturated laces is always the same, one has the inequalities:
$$
N_0(X) \leq  \sum_{L \in \L_s} (-1)^{|L|} N(L) \, \quad or
$$
$$
N_0(X) \geq \sum_{L \in \L_s} (-1)^{|L|} N(L) \,\,\,,
$$
according to whether that parity is odd or even respectively.
 
{\bf EXAMPLES}
 
{\bf 1) The identity lace-map} : $l(S)=S$ for every $S \subset P$.
Every subset $S$ of $P$ is a lace, and $C(S)$ is always empty.
In this case the Lace Expansion reduces to the inclusion-exclusion
principle.
 
{\bf 2) The Bonferroni Lace}:
Let $P=\{1, \dots , n\}$, for a positive integer $n$. 
For $S \subset P$, define
$l(S)$ to be the subset of $S$ consisting
of its $k$ smallest elements, if $|S| \geq k$, and
otherwise $l(S):=S$. It is easy to see that $l$ is
a lace-map, and that the laces are the subsets of $P$ with
cardinality $\leq k$. The saturated laces are those
whose cardinality is $<k$ and for an unsaturated
lace $L=\{i_1 < \dots < i_k\}$, $C(L)={i_k+1, i_k+2 , \dots , n}$.
The resulting sieve is the Bonferroni sieve (see, e.g., [C]).
 
{\bf 3) The Brun Lace} : Let $P= \{ 1, 2, \dots , n \}$,
and let $n \geq N_1 \geq N_2 \geq \dots \geq N_n \geq 1$ be given
beforehand.
For a set $S=\{ i_1 > i_2 > \dots > i_r \}$, define 
$l(S)=S$ if $i_1 \leq N_1 , i_2 \leq N_2 , \dots , i_r \leq N_r$,
and otherwise, let $l(S)= \{ i_1 , \dots , i_s \}$, where $s$ is the
smallest index such that $i_s > N_s$. It is easy to see that $l$ is
a lace-map. The saturated laces are sets of the form
$L=\{ i_1 > \dots > i_r \}$, where
$i_1 \leq N_1 , i_2 \leq N_2 , \dots , i_r \leq N_r$. 
The unsaturated laces are sets $L=\{ i_1 > \dots > i_r \}$, such that
$i_1 \leq N_1 , i_2 \leq N_2 , \dots , i_{r-1} \leq N_{r-1}$, but
$i_r > N_r$. For such a lace, 
$C(L)=\{ i_{r} -1 , i_{r} -2 , \dots , 1 \}$. 
 
To get upper and lower {\it Brun sieves}([B]) we must have
$N_1=N_2, N_3=N_4, \dots N_{n-1}=N_n$, and
$N_2=N_3, N_4=N_5, \dots N_{n-1}=N_n$,  respectively.
 
{\bf 4) The Brydges-Spencer Original Lace Expansion}[BS]
(see [MS] for a very lucid exposition):
Fix $n$, and let $P=\{ (i,j) ; 0 \leq i < j \leq n \}$.
It is instructive to think of $n+1$ dots placed in a row,
at locations $0, 1 , \dots , n$.
Then $P$ is the set of all arcs joining two dots.
For any set of arcs $G$, the lace-map, $l(G)$, is defined as follows.
Let  $s_1 := 0$, and let $t_1$ be the largest $t$ such that
$0t \in G$. Among all the arcs $st \in G$ that ``go over'' $t_1$,
let $t_2$ be the endpoint that goes farthest amongst them:
$t_2 := max \{t : st \in G, s< t_1 \}$, and amongst those
arcs take $s_2$ to be the left-endpoint
farthest to the left that connects to $t_2$:
$s_2 :=min \{ s : s t_2 \in G \}$. One continues recursively:
$t_{i} := max \{t : st \in G, s< t_{i-1} \}$, and 
$s_i :=min \{ s : s t_i \in G \}$, until either one gets,
for some $k$, $t_k=n$, or one 
finishes the {\it connected component} and starts
a new component with the next dot.  The lace of $G$, $l(G)$,
is defined to be the collection of
these arcs $ \{ s_1 t_1 , s_2 t_2 \dots, s_kt_k \}$.
Observe that the arcs of an irreducible lace 
interlace, forming the kind of lace people embroider,
which explains its name.
 
The Brydges-Spencer lace expansion does not yield any sieves, but
serves another purpose: to determine the asymptotic behavior of
the number of
$n$-step self-avoiding walks. Brydges and Spencer[BS] determined the
asymptotics for {\it weakly avoiding} walks in dimensions
$d \geq 5$, while Hara and Slade[HS][MS],
in one of the greatest mathematical feats of this decade,
strengthened it to
the regular self-avoiding walk in $d \geq 5$. The problem is
still open for dimensions $d=2,3,4$. Perhaps we need
a more complicated lace-map.
 
{\bf SPECULATIONS}
 
I predict that the Abstract Lace Expansion (ALE) has a bright future.
The very powerful {\it Probabilistic Method}([ASE],[S]) uses
the Bonferroni sieve with $k=2$. Introducing appropriate lace-maps
may make it even more powerful. The {\it Satisfiability Problem}
(the grandmother of all NP-complete problems) can be approached
via counting (the number of covered
$0$-$1$ vectors), and introducing powerful lace-maps may improve
the average-running-time of current algorithms. Another possible
application is to improving
current asymptotic upper and lower bounds for $R(n,n)$, as well
as to the exact evaluations of $R(5,5)$ and $R(6,6)$, in spite of
Paul Erd\H os's pessimistic prophesy (see [S], p. 4).
Recall that the Ramsey number, $R(n,n)$, is the
smallest $N$ such that if you 2-color the edges of the complete
graph on $N$ vertices, then you are {\it guaranteed} a monochromatic
$K_n$. 
 
This gives rise to the following counting problem.
Let $X$ is the set of all ${{N} \choose {2}}$ edge-colorings, and
let $P$ consist of the ${{N} \choose {n}}$ properties:
$A_S:=$ the induced coloring on $S$ is monochromatic, where $S$
ranges over all $n-$subsets of the set of vertices $\{ 1, \dots , N\}$.
Find the number $N_0(X)$ of property-less colorings.
If, thanks to some lower sieve,
we can ascertain that $N_0(X)>0$, then we would get that $R(n,n)>N$.
If, on the other hand, thanks to an upper sieve, we would find
that $N_0(X) \leq 0$, then we would know that $N_0(X)=0$, and that
$R(n,n) \leq N$. So we need good lace-maps that would give good
sieves that, in turn, would make the lower and upper bounds zero-in at the
exact value.
\eject
{\bf REFERENCES}
 
[ASE] N. Alon and J. Spencer, with P. Erd\H os,
{\it ``The Probabilistic Method''}, Wiley, NY, 1992.
 
[B] V. Brun, {\it Le crible d'Eratosth\`ene et le th\'eor\`eme
de Goldbach}, Skrifter utgit av Videnskapsselskapet i Kristiania, I.
Matematisk-Naturvidenskabelig Klasse {\bf 1}(1920), 1-36.
 
[BS] D.C. Brydges and T. Spencer, {\it Self-avoiding walks
in 5 or more dimensions}, Comm. Math. Phys. {\bf 97}(1985), 125-148.
 
[C] L. Comtet, {\it ``Advanced Combinatorics''},
Dordrecht-Holland/Boston,  1974.
 
[HS] T. Hara and G. Slade, {\it The lace expansion for self-avoiding
walk in five or more dimensions}, Reviews in Math. Phys. {\bf 4}(1992),
235-327.
 
[MS] N. Madras and G. Slade, ``{\it The Self Avoiding Walk}'',
Birkhauser, Boston, 1993.
 
[S] J. Spencer, {\it ``Ten Lectures on the Probabilistic Method''},
SIAM, Philadelphia, 1987.
 
\bye